\documentstyle[amscd, amssymb, 11pt]{amsart}

\newtheorem{theorem}{Theorem}[section]
\newtheorem{lemma}[theorem]{Lemma}

\newtheorem{cor}[theorem]{Corollary}

\newtheorem{conjecture}[theorem]{Conjecture}

\theoremstyle{definition}

\theoremstyle{remark}
\newtheorem{remark}[theorem]{Remark}

\numberwithin{equation}{section}

\newenvironment{prf}
         {\begin{trivlist}\item[]{\bf Proof:\ }}
                  {{\qed} \end{trivlist}}

\renewcommand{\index}{{\rm Index\, }}
\newcommand{\tr}{{\rm tr}}

\begin{document}

\title{On a Proof of the $\frac{10}{8}$-Conjecture}

\author{Jin-Hong Kim}




\address{Department of Mathematics\\
         Oklahoma State University\\
         Stillwater, OK 74078}
\email{jinkim@@math.okstate.edu}

\begin{abstract}
Let $X$ be a smooth closed oriented non-spin 4-manifold with even 
intersection form $kE_8\oplus nH$. In this article we show 
that we should have $n\ge |k|$ on $X$. Thus we confirm the $\frac{10}{8}$-conjecture affirmatively. As an application, we also give an estimate of intersection forms of spin coverings of non-spin 4-manifolds with even intersection forms.
\end{abstract}

\maketitle

\section{Introduction} 

In 4-manifold theory, there are two fundamental questions of the existence and
uniqueness of smooth 4-manifolds with the given intersection forms. 
It turns out that the classification of smooth oriented 4-manifolds up to 
diffeomorphism is very different from the classification of unimodular forms.
In this paper, we investigate which even forms can be realized as 
intersection forms of smooth 4-manifolds. The algebraic classification of 
unimodular indefinite forms is simple. Indeed, the classification of Hasse and Minkowski says that any odd indefinite form is equivalent over the integers to one of the $m\langle 1\rangle\oplus n\langle -1 \rangle$ and any even indefinite form to one of the $mE_8\oplus nH$ ($n\ge 0$), where $E_8$ is the irreducible negative definite even form of rank 8 associated to the Lie algebra $E_8$ and $H$ is the hyperbolic form of rank 2 (see \cite{D-K} or \cite{M-H}). 

In \cite{D}, S.K. Donaldson showed that definite even forms cannot be realized as the intersection forms of smooth 4-manifolds. In other words, among all definite forms, only the diagonalizable ones over the integers can be realized.
In case of the intersection forms of spin 4-manifolds, M. Furuta proved the 
following estimate:

\begin{theorem}[\cite{Br}, \cite{Furuta1}]
Let $X$ be a smooth closed oriented spin 4-manifold with the intersection form 
$2kE_8\oplus nH$. If $n\ne 0$, then we have $n\ge 2|k|+1$.
\end{theorem}

This theorem has been dubbed as the $\frac{10}{8}$-theorem 
for spin 4-manifolds.  The $\frac{10}{8}$-theorem is a partial result to the 
$\frac{11}{8}$-conjecture stating that in the above theorem we should have 
$n\ge 3|k|$.  

Recently R. Lee and T.-J. Li investigated the intersection forms of smooth non-spin 4-manifolds $X$ with even intersection form $kE_8\oplus nH$, and got a partial result to the $\frac{10}{8}$-conjecture below (see \cite{Lee-Li}). Indeed, using the Furuta's $\frac{10}{8}$-theorem and covering tricks, they proved that under the conditions of the torsion part of $H_1(X, {\Bbb Z})$ we should have $n\ge |k|$. Then they made the following conjecture (see also \cite{Bohr}):

\smallskip
\begin{conjecture}[\bf $\frac{10}{8}$-conjecture] \label{conj1.1}
Let $X$ be a smooth closed oriented non-spin 4-manifold with even 
intersection form $kE_8\oplus nH$. Then we have $n\ge |k|$.
\end{conjecture}
\smallskip

The purpose of this paper is to prove the above conjecture completely (see
Theorem \ref{thm3.1}). Note that the Enrique surface is a non-spin 4-manifold with the intersection form $E_8\oplus H$ and that its $n$-fold connected sum is also a non-spin 4-manifold with intersection form $nE_8\oplus nH$. Thus the lower bound of the conjecture is saturated.

One of the main ingredients to prove the conjecture is that we use a 
$Pin(2)$-equivariant map induced from the Seiberg-Witten equations and compute
the $K$-theoretic degree using the idea of J. Bryan in \cite{Br}. In case of
spin 4-manifolds, we can always construct a $Pin(2)$-equivariant map from the 
Seiberg-Witten equations for the trivial $spin^c$ structure. 
On the contrary, in case of non-spin 4-manifolds with even intersection forms, we can only construct an $S^1$-equivariant map, when we fix a $spin^c$ structure.

To circumvent this problem, we use the family of $spin^c$ structures to
construct a $Pin(2)$-equivariant map induced from the Seiberg-Witten 
equations.  In fact, it turns out that it is enough to use only \lq\lq two" $spin^c$ structures in the 2-torsion part of $H^2(X, {\Bbb Z})$ which correspond to each other via the well-known involution in the Seiberg-Witten theory (see \cite{Morgan}).

Once we construct a $Pin(2)$-equivariant map, we can use the finite 
dimensional approximation of Furuta to get a $Pin(2)$-equivariant 
map on the balls preserving the boundaries of the balls.  However, such finite dimensional approximations conatin both the trivial and non-trivial 1-dimensional representation of $Pin(2)$. This prevents us from applying directly the Furuta's argument for spin 4-manifolds to non-spin 4-manifolds with even intersection form. To overcome this difficulty, we use only ${\Bbb Z}_4$-symmetry generated by the above involution and its representations instead of the full $Pin(2)$-symmetry in order to get some non-trivial information on the $K$-theoretic degree induced by the finite dimensional approximation. In other words, applying the tom Dieck's character formula to the map $f^\ast$ on the ${\Bbb Z}_4$-equivariant $K$-theory induced from the ${\Bbb Z}_4$-equivariant map on the balls, we can easily figure out the 
$K$-theoretic degree of $f^\ast$ in the character ring  $R({\Bbb Z}_4)$ of ${\Bbb Z}_4$. Since $R({\Bbb Z}_4)$ is a ring over the integers, all the coefficients of the 
$K$-theoretic degree are integers. This implies our inequality in the 
$\frac{10}{8}$-conjecture.

Combining our result with Donaldson's results on the definite intersection forms and Furuta's $\frac{10}{8}$-theorem, we can give a complete classification of the intersection forms of smooth 4-manifolds as follows:

\begin{theorem} \label{thm1.2}
Let $X$ be a smooth closed oriented 4-manifold with intersection matrix $Q$. Then $Q$ satisfies the following two possibilities:
\begin{itemize}
\item[(1)] If $Q$ is odd, then 
  $Q=m\langle 1 \rangle\oplus n \langle -1 \rangle$ with
  $m, n\ge 0$.
\item[(2)] If $Q$ is even, then $Q=kE_8\oplus nH$ with $n\ge |k|$.
\end{itemize}
\end{theorem}

For the rest of this note, we will assume that $b_1(X)=0$, due to the 
following lemma using a surgery argument:

\begin{lemma} \label{lem1.1}
Let $X$ be a smooth closed oriented 4-manifold with its intersection matrix $Q_X$. Then there exists a smooth closed oriented 4-manifold $X'$ satisfying the following two conditions;
\begin{itemize}
 \item[(1)] $Q_{X'}=Q_X$.
 \item[(2)] $H_1(X', {\Bbb Z})$ consists of only torsion elements.
\end{itemize}
\end{lemma}

\begin{prf} 
See Lemma 2.4 in \cite{Lee-Li} or see \cite{Furuta1}.
\end{prf}
\smallskip

We organize this paper as follows. In Section 2, we set up basic notations,
and prove important facts necessary for the proof of our main theorem \ref{thm3.1}. Section 3 is devoted to proving the main theorem of this paper. Finally we give an application about the intersection forms of spin coverings of
non-spin 4-manifolds with even intersection forms in Section 4.


\section{Monopole Maps for non-spin 4-manifolds} 

In this section, we will fix notations we use in this paper, and prove an important Theorem \ref{thm2.1} necessary for the proof of Theorem \ref{thm3.1}. 

Let ${\Bbb H}$ be the quaternion numbers, $Sp(1)$ the group of quaternions with norm 1, and $S^1$ the intersection of $Sp(1)$ with ${\Bbb C}$ in ${\Bbb H}$.
We then can define five $Spin^c_4$-modules $_-{\Bbb H}_+$, $_+{\Bbb H}$,
$_-{\Bbb H}$, $_+{\Bbb H}_+$, and $\tilde{\Bbb C}$ such that the actions
of $(q_-, q_+, z)\in Spin^c_4=Sp(1)\times Sp(1)\times_{{\Bbb Z}_2}
U(1)$ on $a\in {}_-{\Bbb H}_+$, $\phi\in {}_+{\Bbb H}$, 
$\psi\in {}_-{\Bbb H}$, $\omega\in {}_+{\Bbb H}_+$, and $\eta\in
\tilde{\Bbb C}$ are defined by 
\[
q_-aq_+^{-1}, q_+\phi z^{-1},
q_-\psi z^{-1}, q_+\omega q_+^{-1},
\]
and $\eta z^2$, respectively.
For a $spin^c$ structure $P$ on $M$, we have the five vector bundles $T, 
S_+, S_-, \Lambda$, and $L$ associated to the five $Spin^c_4$-modules 
$_-{\Bbb H}_+$, $_+{\Bbb H}$, $_-{\Bbb H}$, $_+{\Bbb H}_{+}$, 
and ${\tilde {\Bbb 
C}}$, respectively.

Recall that we call the homotopy class of $(P,T)$ a $spin^c$ structure if 
$P\times_{Spin^c_4}{ } _-{\Bbb H}_{+}=T\cong T^\ast X$.  

The $Spin^c_4$-equivariant map $_-{\Bbb H}_{+}\times {}_+{\Bbb H}\to 
{}_-{\Bbb H}$ 
given by $(a,\phi)\mapsto a\phi$ induces Clifford multiplication 
\[
C:T\otimes S_+\to S_-.
\]
Similarly, we can define the twisted Clifford multiplication 
\[
{\tilde C}: T\otimes T\to \Lambda,
\]
using the $Spin^c_4$-equivariant map $_-{\Bbb H}_{+}\times 
{}_-{\Bbb H}_{+}\to 
{}_+{\Bbb H}_{+}$ given by $(a,b)\mapsto \bar a b$. In our case, the twisted 
Clifford multiplication is given by
\[
\bar C:T^\ast X\otimes T^\ast X\to \Lambda^0\oplus \Lambda^+,
(a,b)\mapsto \langle a, b\rangle\oplus p_+(a\wedge b),
\]
where $p_+:\Lambda^2\to \Lambda^+$ is the orthogonal projection and $\langle
\cdot, \cdot \rangle$ is the inner product.

Let $A_0$ be a fixed connection on $L$ associated to the $spin^c$ structure 
$P$. Then we have the twisted Dirac operators
\[
D_1=C\circ \nabla_1: \Gamma(S_+)\to \Gamma(S_-),
\]
and
\[
D_2=d^\ast + d^+=\bar C\circ \nabla_2 :\Gamma(T^\ast X)\to \Gamma(
\Lambda^0)\oplus
\Gamma(\Lambda^+),
\]
where $\nabla_1$ is the covariant derivative on $\Gamma(S_+)$ induced from 
the Riemannian connection on $T^\ast X$ 
and the fixed connection $A_0$ on $L$, and $\nabla_2$ is the covariant
derivative on $\Gamma(T^\ast X)$ induced from the Riemannian connection on
$T^\ast X$.

We define the quadratic map 
\[
Q:\Gamma(S_+)\times \Gamma(T^\ast X)\to \Gamma(S_-)\times\Gamma(\Lambda^+)
\]
induced from the $Spin^c_4$-equivariant map
\[
_+{\Bbb H}\times { }_-{\Bbb H}_{+}\to { }_-{\Bbb H}\times { }_+{\Bbb H}_{+},
(\phi, a)\mapsto (a\phi i, \phi i \bar\phi).
\]
Consider the monopole map
\[
D_{A_0}\oplus Q: V\to W,
\]
where $V$ is the $L^2_4$-completion of $\Gamma(S_+\oplus T^\ast X)$ and $W$ 
is the $L^2_3$-completion of $\Gamma(S_-\oplus \Lambda^0\oplus \Lambda^+)$.
Let $F_{A_0}^+$ be the self-dual part of the curvature of $A_0$. When we
identify the imaginary part of $\Gamma(\Lambda)$ induced from
$_{+}{\Bbb H}_{+}$ with the space of self-dual 2-forms, the 
\emph{monopole equation}
is defined to be $(D+Q)(v)+F_{A_0}^+=0$ for $v\in V$,
and solutions of the monopole equation are called \emph{monopoles}.

The monopole map $D_{A_0}\oplus Q$ has the symmetry of the stabilizer of $A_0$. 
The stabilizer of $A_0$ is the group of the harmonic maps ${\rm Harm}(X, S^1)$
whose group structure is induced from that of $S^1$.
Since we assume that $b_1(X)=0$, the harmonic maps are constants. Thus 
one can see that $D_{A_0}\oplus Q$ is only an 
$S^1$-equivariant map. For the sake of convenience, from now on 
we omit the subscript $A_0$, unless there is some confusion.

Let $V_\lambda^1$ be the subspace of $L_2^4(S_+)$ spanned by the eigenvectors of $D_1^\ast D_1$ with eigenvalues less than or equal to $\lambda$, and let $W_\lambda^1$ be the subspace of $L_2^4(S_-)$ spanned by the eigenvectors of $D_1^\ast D_1$ with eigenvalues less than or equal to $\lambda$.  Similarly, define $V_\lambda^2$ and $W_\lambda^2$ using the elliptic operator $D_2$. 
Let $V_\lambda=V_\lambda^1\oplus V_\lambda^2$ and $W_\lambda=W_\lambda^1\oplus W_\lambda^2$. Using the monopole map $D_{A_0}\oplus Q$, Furuta constructed a finite dimensional approximation, and its finite dimensional approximation 
$D_\lambda\oplus Q_\lambda: V_{\lambda}\to W_{\lambda}$ is also an 
$S^1$-equivariant map satisfying the following two conditions;
(1) for large enough $\lambda$, it does not vanish on the finite dimensional sphere in $V_{\lambda}$ of radius $R$ with center 0, and (2) the image of the sphere is contained in $W_{\lambda}$.

To construct a $Pin(2)$ symmetry, we have to use all or some $spin^c$ 
structures on $X$. Recall that $Pin(2)$ is the normalizer of $S^1$ in 
$Sp(1)$, and it is generated by $S^1$ and $j$. In this note, we use only two 
$spin^c$ structures to construct a $Pin(2)$-equivariant map induced from the 
monopole maps.

To do this, we have to recall the well-known involution map in the 
Seiberg-Witten theory. Let 
\[
\iota: Spin^c_4=Sp(1)\times Sp(1)\times_{{\Bbb Z}_2} U(1)\to Spin^c_4
\]
be the involution given by $(q_-, q_+, z)\mapsto (q_-, q_+, z^{-1})$. Then we 
have the five associated $Spin^c_4$-representations $_-{\Bbb H}_+'$ , 
$_+{\Bbb 
H}'$, $_-{\Bbb H}'$, $_+{\Bbb H}_{+}'$, and ${\tilde {\Bbb C}'}$ obtained by 
twisting by $\iota$  $_-{\Bbb H}_+$, $_+{\Bbb H}$, $_-{\Bbb H}$, 
$_+{\Bbb H}_{+}$, 
and ${\tilde {\Bbb C}}$, respectively.  Let $P'=P\times_{\iota} Spin^c_4$ be 
the twisting of $P$ by $\iota$, and $T', S_+', S_-', \Lambda'$, and $L'$ denote 
the five associated vector bundles to the five $Spin^c_4$-modules $_-{\Bbb 
H}_+'$, $_+{\Bbb H}'$, $_-{\Bbb H}'$, $_+{\Bbb H}_{+}'$, 
and ${\tilde {\Bbb C}}'$, 
respectively.
Note that the five $Spin^c_4$-modules $_-{\Bbb H}_+$, $_+{\Bbb H}$, 
$_-{\Bbb H}$, 
$_+{\Bbb H}_{+}$, and ${\tilde {\Bbb C}}$ correspond to the five 
$Spin^c_4$-modules $_-{\Bbb H}_+'$, $_+{\Bbb H}'$, $_-{\Bbb H}'$, $_+{\Bbb 
H}_{+}'$, and ${\tilde {\Bbb C}}'$ via the following $Spin^c_4$-equivariant 
homomorphisms
\begin{equation} \label{eq2.3}
\begin{split}
 _-j_{+} &: {}_-{\Bbb H}_+\to {}_-{\Bbb H}_+',\quad  a\mapsto -a\\
 _+j &: {}_+{\Bbb H}\to {}_+{\Bbb H}',\quad  \phi\mapsto \phi j\\
 _-j &: {}_-{\Bbb H}\to {}_-{\Bbb H}',\quad  \psi\mapsto \psi j\\
 _+j_{+} &: {}_+{\Bbb H}_+\to {}_+{\Bbb H}_+',\quad  \omega\mapsto -\omega\\
 \tilde j &: \tilde {\Bbb C}\to \tilde {\Bbb C}', \quad t\mapsto \bar t.\\
\end{split}
\end{equation}

Let $Spin^c_4(X)$ denote the set of all inequivalent $spin^c$ structures on 
$X$. Note that the twisting of $\iota$ induces an involution on $Spin^c_4(X)$. In particular, the twisting of $\iota$ induces an involution on the set $\{ L, 
L'=L^{-1} \}$. Moreover, if $L$ is the trivial line bundle which is the case of the spin structure or $L$ is in the 2-torsion part of $H^2(X, {\Bbb Z})$, then the twisting of $\iota$ induces a trivial involution on the set $\{ L, L' \}$ (e.g., see Lemma 3.7 in \cite{Furuta2}).

From now on, we will use the superscript prime to denote any objects of $P'$  
corresponding to $P$. Let $A_0$ and $A_0'$ be connections on $L$ and $L'$ which corresponds to each other by the map $\tilde j$. Let $J$ denote the bundle maps induced from the five isomorphisms in \eqref{eq2.3}. Then from the above discussion we have the following:

\begin{lemma} \label{lem2.1}
The two maps $D\oplus Q:V\to W$ and $D'\oplus Q':V'\to W'$ correspond to each 
other via the map $J$ as follows:
\begin{equation*} 
\begin{CD}
V @>\quad\quad J\quad\quad >> V'\\
@V{D_{A_0}\oplus Q}VV   @VV{D'_{A'_0}\oplus Q'}V\\
W @>\quad\quad J \quad\quad >> W'.\\
\end{CD}
\end{equation*}
\end{lemma}

In case a $spin^c$ structure is reduced to the spin structure, we take the 
trivial product connection $A_0$, and the map $D\oplus Q$ and its finite 
dimensional approximation have a $Pin(2)$ symmetry. Thus Furuta was able to 
prove the $\frac{10}{8}$-conjecture on spin 4-manifolds (see \cite{Furuta1}). 
Later, Bryan reproved and improved his results using the same $Pin(2)$ 
symmetry and the tom Dieck's character formula instead of the Adams 
operations, which is the method of this paper (see \cite{Br} for more 
details).  Note also that even if the twisting of $\iota$ induces a trivial involution on the set $\{ L, L' \}$ for a 2-torsion class $L$, the map $D\oplus Q$ do not admits a $Pin(2)$-symmetry in general, since we cannot take the product trivial connection as a base connection $A_0$.

We will also use the following lemma in the proof of the Theorem 
\ref{thm2.1}.

\begin{lemma} \label{lem2.2}
Let $X$ be a smooth oriented non-spin 4-manifold with even intersection 
form. Then $w_2(X)$ has an integral lift which is in the 2-torsion part of $H^2(X,{\Bbb Z})$.
\end{lemma}

\begin{remark}
In particular, this lemma implies that $H_1(X, {\Bbb Z})$ must have ${\Bbb Z}_{2^i}$ as a direct summand of $H_1(X, {\Bbb Z})$.
\end{remark}

\begin{prf} This is the Lemma 2.2 in \cite{Lee-Li} (see also \cite{A-L}). For the sake of completeness, we give its proof.

We first consider the following universal coefficient sequences with integral
and ${\Bbb Z}_2$ coefficients
\begin{equation*} 
\begin{CD}
\text{Ext}(H_1(X, {\Bbb Z}), {\Bbb Z})@>k >> H^2(X, {\Bbb Z}) 
@>h_1 >>\text{Hom}(H_2(X, {\Bbb Z}), {\Bbb Z})\\
@V{\rho}VV   @V VV   @V VV  \\
\text{Ext}(H_1(X, {\Bbb Z}), {\Bbb Z}_2)@>  >>H^2(X, {\Bbb Z}_2)
@> h_2 >>\text{Hom}(H_2(X, {\Bbb Z}), {\Bbb Z}_2).\\
\end{CD}
\end{equation*}
Recall that the homomorphisms $h_1$ and $h_2$ are related to the intersection
form
\begin{equation}\label{eq2.2}
h_1(\alpha)(b)=a\cdot b, \quad h_2(\alpha)(b)\equiv a\cdot b
\quad {\rm mod}\ 2,
\end{equation}
where $\alpha$ is the Poincar\' e dual of $a$ with either ${\Bbb Z}$ or
${\Bbb Z}_2$ coefficients.
Since the intersection form of $X$ is even, it follows from \eqref{eq2.2} that
$h_2(w_2(X))$ is trivial. Thus $w_2(X)$ comes from a unique element $u\in
\text{Ext}(H_1(X, {\Bbb Z}), {\Bbb Z}_2)$. Since ${\Bbb Z}$ is mapped onto 
${\Bbb Z}_2$, the map $\rho$ is also onto. Thus we can choose $v$ in the
2-torsion part of $\text{Ext}(H_1(X, {\Bbb Z}), {\Bbb Z})$ such that 
$\rho(v)=u$. Let $x=k(v)$. Then $x$ is an integral lift of $w_2(X)$ and is in 
the 2-torsion part of $H^2(X, {\Bbb Z})$ as in the lemma. This completes the proof.
\end{prf}

As a consequence, one can choose a $spin^c$ structure $L$ such that $c_1(L)$ 
is in the torsion part of $H^2(X,{\Bbb Z})$. Note also that $c_1(L)^2$ vanishes.

Let ${\tilde{\Bbb R}}$ be the unique non-trivial 1-dimensional 
representation of $Pin(2)$, and ${\Bbb H}$ the representation which is the restriction of the standard representation of $SU(2)=Sp(1)$ to $Pin(2)\subset Sp(1)$.

The main purpose of this section is to prove the following theorem:

\begin{theorem} \label{thm2.1}
Let $X$ be a smooth closed oriented non-spin 4-manifold with even intersection form and $b_1(X)=0$. Assume that the signature $\sigma(X)$ of $X$ is non-negative and  that $k=-\frac{\sigma(X)}{8}$. Then there are a finite dimensional real $Pin(2)$-modules ${\cal V}_\lambda$ and ${\cal W}_\lambda$, a 
$Pin(2)$-equivariant linear map ${\cal D}_\lambda$, and a 
$Pin(2)$-equivariant quadratic map ${\cal Q}_\lambda$ from ${\cal V}_\lambda$ to ${\cal W}_\lambda$ which satisfy the following conditions:
\begin{itemize}
 \item[(1)] There are $Pin(2)$-module isomorphisms
${\cal V}_\lambda={\Bbb H}^{k+m}\oplus {\tilde{\Bbb R}}^{n}\oplus {\Bbb R}^{n}$ and $\bar{\cal W}_\lambda={\Bbb H}^m\oplus {\tilde{\Bbb R}}^{b_+ +1+ n}\oplus {\Bbb R}^{b_+ +n}$ for some non-negative $m$ and $n$.
 
 \item[(2)] For large enough $\lambda$, ${\cal D}_\lambda+ {\cal Q}_\lambda$ 
does not vanish on the finite dimensional sphere in ${\cal V}_{\lambda}$ of 
radius $R$ with center 0 defined by using a $Pin(2)$-invariant metric on
${\cal V}_\lambda$, and the image of the sphere is contained in 
$\bar{\cal W}_{\lambda}$.
\end{itemize}
\end{theorem}

\begin{prf}
We first take a $spin^c$ structure $L$ such that $c_1(L)$ is in the 2-torsion
part of $H^2(X, {\Bbb Z})$ by Lemma \ref{lem2.2}. The choice of such a 
$spin^c$ structure has two important consequences: first, $c_1(L)^2$ vanishes. Secondly,  the monopole equation for such $spin^c$ structures has a flat connection and zero spinor as a reducible solution so that we can construct the finite dimensional approximation. Note also that since $L$ is in the 2-torsion part of $H^2(X, {\Bbb Z})$, $L=L'$. Nonetheless, we need two copies of monopole maps to construct a $Pin(2)$-equivariant map, since we need two connections $A_0$ and $A_0'$ on $L$ and $L'=L$ which correspond to each other via the map $\tilde j$.

To begin with our proof, we consider the map 
\[
{\cal D}+{\cal Q}: {\cal V}\to {\cal W}, ((\phi, a), (\psi, b))\mapsto
 ((D_{A_0}\oplus Q)(\phi, a), (D'_{A_0'}\oplus Q')(\psi, b)),
\]
where
\begin{equation*}
\begin{split}
{\cal V} &:= \Gamma(S_+\oplus T^\ast X)\times \Gamma(S_+'\oplus T^\ast X),\\
{\cal W} &:= \Gamma(S_-\oplus \Lambda^0(X)\oplus \Lambda^+(X))\times
             \Gamma(S_-'\oplus \Lambda^0(X)\oplus \Lambda^+(X)). \\
\end{split}
\end{equation*}
Using the finite dimensional approximations of $D\oplus Q$ and $D'\oplus Q'$, 
we have a good finite dimensional approximation ${\cal D}_\lambda+{\cal 
Q}_\lambda: {\cal V}_\lambda\to {\cal W}_\lambda$ in that for large enough 
$\lambda$, it does not vanish on the sphere $S_R(V_\lambda\oplus V_\lambda')$ 
of radius $R$ with center 0 and  its image of the sphere is contained in the 
finite dimensional vector space $W_\lambda\oplus W'_\lambda$. 
 
Note from Lemma \ref{lem2.1} that we can define a self map $J$ on  
${\cal V}$ and ${\cal W}$ induced from the five maps in \eqref{eq2.3}, and thus on their finite dimensional approximations so that ${\cal D}+{\cal Q}$ and its finite dimensional approximations are $Pin(2)$-equivariant. This action $J$ is not an involution, but is of order 4. Combining the obvious $S^1$ actions on ${\cal V}$ and ${\cal W}$ with the action $J$ induces a $Pin(2)$ symmetry on ${\cal V}$ and ${\cal W}$. It is clear that ${\cal D}+{\cal Q}$ and its finite dimensional approximation are also $Pin(2)$-equivariant maps.

Let ${\cal D}_1= D_1\oplus D_1'$, and let ${\cal V}^1_\lambda$ be the 
subspace of $L_2^4(S_+)\oplus L_2^4(S_+')$ spanned by the eigenspaces of 
${\cal D}_1^\ast{\cal D}_1$ less than or equal to $\lambda$. Similarly, 
define ${\cal D}_2, {\cal W}^1_\lambda, {\cal V}^2_\lambda$, and ${\cal 
W}^2_\lambda$. Then we have the following decompositions as $Pin(2)$-modules
\begin{equation*}
{\cal V}_\lambda={\cal V}^1_\lambda\oplus {\cal V}^2_\lambda, \quad
{\cal W}_\lambda={\cal W}^1_\lambda\oplus {\cal W}^2_\lambda,
\end{equation*}
where ${\cal V}^1_\lambda=V^1_\lambda\oplus {V^1_\lambda}'$, and so on.

We now show that there exists non-negative integers $m$ and $n$ such that as $Pin(2)$-modules
\[
{\cal V}_\lambda={\Bbb H}^{m+k}\oplus {\tilde{\Bbb R}}^{n}\oplus {\Bbb R}^n
\quad\text{and}\quad {\cal W}_\lambda={\Bbb H}^m\oplus {\tilde{\Bbb R}}^{n+1+b_+}\oplus {\Bbb R}^{n+1+b_+}.
\]
To do this, let ${\Bbb C}$ denote the standard complex 1-dimensional representation space of $S^1$. We can choose finitely many points $p_1,p_2, \ldots, p_l$ such that the restriction on fibers over these points is an injection form $V^1_\lambda$ to $\oplus^{l}_{j=1} (S_+)_{p_j}$ which is isomorphic to ${\Bbb H}^l$ as a 
$Pin(2)$-module. 
Since $V_\lambda^1$ is just an $S^1$-module, $V^1_\lambda$ is isomorphic to ${\Bbb C}^{m'}$ for some $m'$ as an $S^1$-module. 
Similarly we may assume that over the same finitely many points $p_1, p_2, \ldots, p_l$ the restriction on fibers induces an isomorphism between ${V^1_\lambda}'$ and ${\Bbb C}^{m'}$ as an $S^1$-module.

Using the map induced by the identification $(\phi, \psi)\mapsto z+wj$ for a pair of generators $\phi\in V^1_\lambda$ and $\psi\in {V^1_\lambda}'$, we next identify ${\cal V}^1_\lambda$ with $({\Bbb C}\oplus {\Bbb C})^{m'}={\Bbb H}^{m'}\subset ({\Bbb H}\oplus {\Bbb H})^{m'}$. Then the map $J$ on ${\cal V}^1_\lambda$ gives rise to a self map, denoted by the same letter $J$, on ${\Bbb H}^{m'}$ induced by $J: z+wj\mapsto -w+zj$ on each component ${\Bbb H}$.  It is also easy to see that $S^1$ acts on ${\Bbb H}^{m'}$ from the right. Thus we can identify ${\cal V}^1_\lambda$ with ${\Bbb H}^{m'}$ 
as a $Pin(2)$-moidule. Similarly we identify ${\cal W}^1_\lambda$ with ${\Bbb H}^m$ with the standard $Pin(2)$-symmetry. 
But, the index of ${\cal D}_1$ equals 
\begin{equation*}
\begin{split}
\dim {\cal V}^1_\lambda-\dim {\cal W}^1_\lambda &=4m'-4m\\
&=\index D_1 +\index D_1'\\
&=\frac{c_1(L)^2-\sigma(X)}{4}+\frac{c_1(L')^2-\sigma(X)}{4}=4k,
\end{split}
\end{equation*}
where we used $c_1(L)^2=c_1(L')^2=0$.
Thus we have $m'=m+k$.

On the other hand, since $J$ acts on ${\cal V}^2_\lambda$ as a map given by $(a, b)\mapsto (-b, -a)$, $J$ has the eigenvalues $\pm 1$. Thus ${\cal V}^2_\lambda$ contains both the trivial representation ${\Bbb R}$ and the non-trivial representation $\tilde{\Bbb R}$ of $Pin(2)$ as summands.  As above, we can see that ${\cal V}^2_\lambda$ is of the form $\tilde{\Bbb R}^n\oplus {\Bbb R}^n$ 
as $Pin(2)$-modules for some $n$. Similarly ${\cal W}^2_\lambda$ is of the form  ${\tilde{\Bbb R}^{n'}}\oplus {\Bbb R}^{n'}$ as $Pin(2)$-modules for some $n'$. But, the index of ${\cal D}_2$ equals 
\[
\index D_2 +\index D_2'=-2(1+b_+)=2n-2n'.
\]
Thus we have $n'=n+(1+b_+)$.  

As in \cite{Furuta1}, we next show that the image of ${\cal D}+
{\cal Q}$ is contained in a subspace of ${\cal W}$ of codimension 1. 
To do this, let $s_0$ be the parallel section of $\Gamma(\Lambda^0\oplus \Lambda^+)$ in $\Gamma(S_-\oplus \Lambda^0 \oplus \Lambda^+)$ which corresponds to the $Spin^c_4$-invariant element $1\in {}_+{\Bbb H}_+$. Since the 
$Spin^c_4$-invariant element $1\in {}_+{\Bbb H}_+$ corresponds to the 
$Spin^c_4$-invariant element $-1\in {}_+{\Bbb H}_+'$ via the map ${}_+ j_+$, we should have the parallel section $s_0'$ of $\Gamma(\Lambda^0\oplus \Lambda^+)$ in $\Gamma(S_-'\oplus \Lambda^0 \oplus \Lambda^+)$ which corresponds to the $Spin^c_4$-invariant element $-1\in {}_+{\Bbb H}_+'$. Thus we have the parallel section $(s_0, s_0')$ of $\Gamma(\Lambda^0\oplus \Lambda^+)\times \Gamma(\Lambda^0\oplus \Lambda^+)$ in $\Gamma(S_-\oplus \Lambda^0 \oplus \Lambda^+)\times\Gamma(S_-'\oplus \Lambda^0 \oplus \Lambda^+)$, and it corresponds to the $Spin^c_4$-invariant element $(1, -1)\in {}_+{\Bbb H}_+\oplus {}_+{\Bbb H}_+'$.

Since parallel sections are in the kernel of the 
${\cal D}^\ast$ and the image of ${\cal D}$ is $L^2$-orthogonal to the kernel 
of ${\cal D}^\ast$, the image of ${\cal D}$ should be contained in the 
$L^2$-orthogonal complement $(s_0, s_0')^\perp$ of $(s_0, s_0')$ in the space $W$. By the construction of ${\cal Q}$, it is also clear that the image of ${\cal Q}$ is contained in the $(s_0, s_0')^\perp$.  

To finish the proof of this theorem, note from the construction of the finite 
dimensional approximation that the image of its finite dimensional 
approximation ${\cal D}_\lambda+{\cal Q}_\lambda$ is also contained in a 
subspace of $\bar{\cal W}_\lambda={\cal W}_\lambda\cap (s_0, s_0')^\perp$ of 
codimension 1. Since ${\Bbb R}(s_0, s_0')$ corresponds to the 
$Spin^c_4$-invariant element $(1, -1)\in {}_+{\Bbb H}_+\oplus {}_+{\Bbb H}_+'$, it is the trivial representation ${\Bbb R}$ as a $Pin(2)$-module. Thus
we have $\bar{\cal W}_\lambda={\Bbb H}^m\oplus \tilde{\Bbb R}^{b_++1+n}\oplus{\Bbb R}^{b_+ + n}$. This completes the proof of (1) and (2).
\end{prf}

As we remarked in the introduction, we cannot extract any non-trivial information about the intersection forms of non-spin even 4-manifolds from Theorem \ref{thm2.1}, because of the 1-dimensional trivial representations of $Pin(2)$.  To overcome this difficulty by removing the 
1-dimensional trivial representations of $Pin(2)$ in the above theorem, we use only ${\Bbb Z}_4$-symmetry induced by the map $J$ as follows.

\begin{cor} \label{cor2.1}
Let $X$ and $k$ be the same as in Theorem \ref{thm2.1}. Then there are a finite dimensional real ${\Bbb Z}_4$-modules ${\Bbb V}_\lambda$ and ${\Bbb W}_\lambda$, a ${\Bbb Z}_4$-equivariant linear map ${\Bbb D}_\lambda$, and a 
${\Bbb Z}_4$-equivariant quadratic map ${\Bbb Q}_\lambda$ from ${\Bbb V}_\lambda$ to ${\Bbb W}_\lambda$ which satisfy the following conditions:
\begin{itemize}
 \item[(1)] There are ${\Bbb Z}_4$-module isomorphisms
${\Bbb V}_\lambda={\Bbb C}^{k+m}\oplus {\tilde{\Bbb R}}^{n}$ and 
${\Bbb W}_\lambda={\Bbb C}^m\oplus {\tilde{\Bbb R}}^{b_+ +1+ n}$ for some 
non-negative $m$ and $n$, where the induced map $J$ on ${\Bbb C}^{m}$ and ${\Bbb C}^{m+k}$ acts as a left multiplication by $i$.
 
 \item[(2)] For large enough $\lambda$, ${\Bbb D}_\lambda+ {\Bbb Q}_\lambda$ 
does not vanish on the finite dimensional sphere in ${\Bbb V}_{\lambda}$ of 
radius $R$ with center 0 defined by using a ${\Bbb Z}_4$-invariant metric on
${\Bbb V}_\lambda$, and the image of the sphere is contained in 
${\Bbb W}_{\lambda}$. 
\end{itemize}
\end{cor}

\begin{prf}
We continue to use the notations in Theorem \ref{thm2.1}.
Let ${\Bbb V}^1_\lambda$ be the $+i$-eigenspace of ${\cal V}^1_\lambda$ of the map $J$.  Let ${\Bbb V}_\lambda={\Bbb V}^1_\lambda\oplus \tilde{\Bbb R}^{n}$ be the subspace of ${\Bbb H}^{m+k}\oplus \tilde{\Bbb R}^{n}$ in ${\cal V}_\lambda$ under the identification. Similarly, let ${\Bbb W}^1_\lambda$ be the 
$+i$-eigenspace of $\bar{\cal W}^1_\lambda$ of the map $J$. Let ${\Bbb W}_\lambda={\Bbb W}^1_\lambda\oplus \tilde{\Bbb R}^{n+1+b_+}$. Then ${\Bbb V}_\lambda$ and ${\Bbb W}_\lambda$ are ${\Bbb Z}_4$-modules of the cyclic group of order 4 generated by the map $J$. It is also clear that ${\Bbb V}^1_\lambda$ and ${\Bbb W}^1_\lambda$ are isomorphic to ${\Bbb C}^{m+k}$ and ${\Bbb C}^{m}$ as ${\Bbb Z}_4$-modules in ${\Bbb H}^{m+k}$ and ${\Bbb H}^{m}$, respectively, where the induced map $J$ acts as a left multiplication by $i$. 

Let ${\Bbb D}_\lambda$ and ${\Bbb Q}_\lambda$ be the restrictions 
of ${\cal D}_\lambda$ and ${\cal Q}_\lambda$ to the 
subspace ${\Bbb V}_\lambda$ of ${\cal V}_\lambda$, respectively. 
Then we can see that the image of the map 
${\Bbb D}_\lambda\oplus {\Bbb Q}_\lambda$ is contained in ${\Bbb W}_\lambda$. 
Finally, since the map ${\cal D}_\lambda\oplus {\cal Q}_\lambda$ does not vanish on the sphere $S_R({\cal V_\lambda})$, the restriction on the subspace does not vanish on the sphere $S_R({\Bbb V}_\lambda)={\Bbb V}_\lambda\cap S_R({\cal V}_\lambda)$, either.  This completes the proof.
\end{prf}

\section{Proof of the $\frac{10}{8}$-Conjecture}  

The purpose of this section is to prove the main theorem. To do it, we use the tom Dieck's character formula for the $K$-theoretic degree. We will briefly explain it (see \cite{Br} or \cite{Dieck} for more details). 

Let $V$ and $W$ be complex $G$ representations for some compact Lie group $G$. Let $BV$ and $BW$ denote balls in $V$ and $W$, and let $f: BV\to BW$ be a $G$-equivariant map preserving the boundaries $SV$ and $SW$. Let $K_G(V)$ denote $K_G(BV, SV)$. Simliarly define $K_G(W)$. Then $K_G(V)$ (resp. $K_G(W)$) is a free $R(G)$-modules with one generator $\lambda(V)$ (resp. $\lambda(W)$), called the \emph{Bott class}.  Now, applying $K$-theory functor to $f$ we get a map
\[
f^\ast: K_G(W)\to K_G(V)
\]
with a unique element $\alpha_f$, called the \emph{$K$-theoretic degree of $f$}, satisfying the equation 
\[
f^\ast(\lambda(W))=\alpha_f\cdot \lambda(V).
\]
Let $V_g$ and $W_g$ denote the subspaces of $V$ and $W$ fixed by an element $g\in G$, and $V_g^\perp$ and $W_g^\perp$ denote the their orthogonal complements. Let $f^g: V_g\to W_g$ be the restriction of $f$ to $V_g$, and let $d(f^g)$ denote the topological degree of $f^g$. Then the tom Dieck's charcter formula says that we have
\begin{equation}\label{eq3.1}
\tr_g(\alpha_f)=d(f^g)\tr_g(\sum_{i=0}^\infty (-1)^i
\Lambda^i(W_g^\perp-V_g^\perp)),
\end{equation}
where $\tr_g$ is the trace of the action of an element $g\in G$.
Note also that the topological degree $d(f^g)$ is by definition zero, if $\dim(V_g)\ne \dim(W_g)$.

\begin{theorem} \label{thm3.1}
Let $X$ be a smooth closed oriented non-spin 4-manifold with even intersection form $kE_8\oplus nH$. Then we have $n\ge |k|$.
\end{theorem}

\begin{prf} 
We assume without loss of generality that the signature $\sigma(X)$ of $X$ is 
non-positive. Note that the complexifications of the spaces ${\Bbb V}_\lambda$ and ${\Bbb W}_\lambda$ are as ${\Bbb Z}_4$-modules
\begin{equation*}
{\Bbb V}_\lambda\otimes_{\Bbb R} {\Bbb C}={\Bbb C}^{m+k}\oplus \bar{{\Bbb C}}^{m+k}\oplus \tilde{{\Bbb C}}^{n}\quad\text{and}\quad
{\Bbb W}_\lambda\otimes_{\Bbb R} {\Bbb C}={\Bbb C}^{m}\oplus \bar{{\Bbb C}}^{m}\oplus \tilde{{\Bbb C}}^{n+1+b_+},
\end{equation*}
where $J$ acts on $\bar{{\Bbb C}}^m$ and $\bar{{\Bbb C}}^{m+k}$ as a left multiplication by $-i$.

Using the finite dimensional approximations 
${\Bbb D}_\lambda+{\Bbb Q}_\lambda$ of Corollary \ref{cor2.1}, we can also construct a ${\Bbb Z}_4$-equivariant map
\[
f_\lambda: (B({\Bbb V}_{\lambda}\otimes {\Bbb C}), S({\Bbb V}_\lambda\otimes 
{\Bbb C}))\to (B({\Bbb W}_\lambda\otimes {\Bbb C}), S({\Bbb 
W}_\lambda\otimes{\Bbb C})).
\]

For the sake of our convenience, we drop the subscript $\lambda$ from now on.

Let ${\Bbb V}_{\Bbb C}={\Bbb V}\otimes {\Bbb C}$ and 
${\Bbb W}_{\Bbb C}={\Bbb W}\otimes {\Bbb C}$.
Applying the $K$-theory functor to $f=f_\lambda$, we get 
\[
f^\ast : K_{{\Bbb Z}_4}({\Bbb W}_{\Bbb C})\to K_{{\Bbb Z}_4}({\Bbb V}_{\Bbb C}),
\]
and there exists a unique element $\alpha_f\in R({\Bbb Z}_4)$ such that 
$f^\ast(\lambda({\Bbb W}_{\Bbb C}))=\alpha_f\cdot \lambda({\Bbb V}_{\Bbb 
C})$.

Since $\alpha_f\in R({\Bbb Z}_4)={\Bbb Z}[\xi]/(\xi^4=1)$, it has the form
\[
\alpha_f=a_01+a_1 \xi+a_2 \xi^2+a_3 \xi^3,\quad a_0,\ldots, a_3\in {\Bbb Z}.
\]
Since $J^2$ acts non-trivially on ${\Bbb C}$, and acts trivially on $\tilde{\Bbb C}$, clearly $\dim ({\Bbb V}_{\Bbb C})_{J^2}\ne \dim ({\Bbb 
W}_{\Bbb C})_{J^2}$. Thus the topological degree $d(f^{J^2})$ is zero.
Applying the tom Dieck's character formula \eqref{eq3.1}, we get
\begin{equation} \label{eq3.2}
0=\tr_{J^2}(\alpha_f)=a_0-a_1 +a_2-a_3.
\end{equation}

On the other hand, since $J$ acts non-trivially on both ${\Bbb C}$ and 
$\tilde{\Bbb C}$, we have $\dim ({\Bbb V}_{\Bbb C})_J=\dim (\bar{\Bbb W}_{\Bbb 
C})_J=0$, and thus the topological degree $d(f^J)$ is 1. By the character formula again we have
\begin{equation*}
\begin{split}
(a_0-a_2)+(a_1-a_3)i &=\tr_J(\alpha_f)=\tr_J(\sum_{i=0}^{\infty}(-1)^i \Lambda^i 
((1+b_+){\tilde{\Bbb C}}-k{\Bbb C}-k\bar{\Bbb C}))\\
&=\tr_J ((1-{\tilde{\Bbb C}})^{1+b_+}(1-{\Bbb C})^{-k} (1-\bar{\Bbb C})^{-k})\\
&=2^{1+b_+} (1-i)^{-k} (1+i)^{-k}\\
&=2^{1+b_+ - k},\\
\end{split}
\end{equation*}
which implies that we have $a_1=a_3$ and $a_0-a_2=2^{1+b_+ -k}$.  Using the 
equation \eqref{eq3.2}, it is easy to see that  $a_0+a_2=a_1+a_3=2a_1$. Thus we have 
\[
2^{1+b_+ -k}=a_0-a_2=2(a_1- a_2).
\]
Hence, $a_1-a_2=2^{b_+-k}$. Since $a_1-a_2$ is an integer, 
we should have $n=b_+\ge k$.  This completes the proof of the $\frac{10}{8}$-conjecture.
\end{prf}

\section{Application: Intersection forms of Spin Coverings} 

The purpose of this section is to apply our main Theorem \ref{thm3.1} to the 
intersection forms of spin covering 4-manifolds of non-spin 4-manifolds with
even intersection forms. 

Recently R. Lee and T.-J. Li proved that every smooth non-spin 4-manifold with even intersection form has a $2^p$-fold spin covering (see \cite{Lee-Li} or \cite{Bohr} for more details). Thus, using the multiplicative property of the signature and Euler characteristic we can prove the following estimate about the intersection forms of spin covering 4-manifolds.

\begin{theorem} \label{thm4.1}
Let $X$ be a smooth closed oriented non-spin 4-manifold with even intersection form, and $M$ its $2^p$-fold covering with intersection form $2kE_8\oplus nH$. Then we have $n\ge 2|k|+2^p-1+b_1(M)-2^p b_1(X)$.
\end{theorem}

\begin{prf}
We assume that
$X$ has the intersection form $lE_8\oplus mH$ with $l,m \ge 0$.
Then the signature and Euler characteristic of $X$ are given by
\[
\sigma(X)=-8l, \quad \chi(X)=2-2b_1(X)+8l+2m.
\]
Using the multiplicative
property of the signature and Euler characteristic between $X$ and $M$, we 
have
\begin{equation} \label{eq4.1}
\sigma(M)=2^p \sigma(X), \quad \chi(M)=2^p \chi(X).
\end{equation}
Since $M$ has the intersection form $2kE_8\oplus nH$, 
it is easy to see from \eqref{eq4.1} 
that $2k=2^p l$ and $n=2^p m +2^p-1+b_1(M)-2^p b_1(X)$.
Now using the Theorem \ref{thm3.1}, it is immediate to get the desired
inequality of the theorem.
\end{prf}

Now we close this section with a remark: if $b_1(X)=b_1(M)=0$, then we have 
the inequality $n\ge 2|k|+2^p-1$. For example, if $H_1(X, {\Bbb Z})$ is 
isomorphic to ${\Bbb Z}_{2^p}$ and the spin covering corresponds to the 
commutator subgroup of $\pi_1(X)$, then it is known that $H_1(M, {\Bbb Z})$ 
consists of only odd torsion elements. Thus we have $b_1(M)=b_1(X)=0$ (see 
\cite{M}). Note that in this case the generator of the deck transformations 
is a spin action of odd type, since it acts freely on $M$ and the quotient 
manifold $X$ is not spin. It is also easy to show that the non-degeneracy 
condition $b_+(X_i)\ne b_+(X_{i-1})$ ($i=2,3, \ldots, p$)
of Theorem 1.2 in \cite{Br} is automatically satisfied. Thus, if either 
$b_+(X)>1$ or $p\ge 4$,  then we can get the inequality  $n\ge
2|k|+p+1$ of Bryan which has been the best known estimate until now.
Hence in certain cases the inequality of Theorem \ref{thm4.1} seems to be 
quite a good improvement.

\end{document}